\documentclass[a4paper,reqno,12pt]{amsart}

\usepackage{enumitem}
\usepackage{layouts}
\usepackage{amsmath, amsthm, amssymb}
\usepackage{amsfonts,mathtools}
\usepackage{bbm} 
\usepackage{epsfig,mathrsfs,graphics}
\usepackage{dsfont,graphicx,adjustbox,upgreek,wrapfig,caption,subcaption,tikz}
\usepackage{epstopdf,latexsym,colortbl,shuffle}

\usepackage[all]{xy} \entrymodifiers={!!<0pt,0.7ex>+}

\usepackage{datetime2}
\usepackage{url}

\usepackage[hyperindex=true, 
					hidelinks]{hyperref}

\newtheorem{thm}{Theorem}
\newtheorem{cor}[thm]{Corollary}

\theoremstyle{definition}
\newtheorem{rem}{Remark}[section]

\numberwithin{equation}{section}

\newcounter{parag}[subsection]

\newcounter{parage}[section]

\newcounter{paraga}

\newcommand{\med}{\medskip}
\newcommand{\bg}{\bigskip}
\newcommand{\sm}{\smallskip}

\newcommand{\al}{\alpha}
\newcommand{\be}{\beta}
\newcommand{\ga}{\gamma}
\newcommand{\Ga}{{\Gamma}}
\newcommand{\de}{\delta}

\newcommand{\De}{\Delta}
\newcommand{\eps}{{\varepsilon}}

\newcommand{\la}{\lambda}
\newcommand{\La}{\Lambda}
\newcommand{\om}{\omega}

\newcommand{\sig}{{\sigma}}

\newcommand{\tht}{{\theta}}
\newcommand{\Th}{\Theta}
\newcommand{\ph}{\varphi}

\newcommand{\demi}{\frac{1}{2}}
\newcommand{\dem}{\tfrac{1}{2}}
\newcommand{\tdemi}{\frac{3}{2}}

\newcommand{\fdemi}{\frac{5}{2}}
\newcommand{\sdemi}{\frac{7}{2}}
\newcommand{\ndemi}{\frac{9
  }{2}}

\newcommand{\ov}{\overline}

\newcommand{\I}{{\mathrm i}}
\newcommand{\dd}{{\mathrm d}}

\newcommand{\ee}{\mathrm e}

\newcommand{\pa}{\partial}

\newcommand{\ii}{^{-1}}

\newcommand{\IM}{\mathop{\Im m}\nolimits}
\newcommand{\RE}{\mathop{\Re e}\nolimits}

\newcommand{\modfour}{\;\textrm{mod}\,4}
\newcommand{\modM}{\;\textrm{mod}\,M}
\newcommand{\modMal}{\;\textrm{mod}\,M_\al}
\newcommand{\modMmual}{\;\textrm{mod}\,M_{\frac{-1}{\al}}}
\newcommand{\modcM}{\;\textrm{mod}\,cM}
\newcommand{\moddM}{\;\textrm{mod}\,2M}

\newcommand{\ie}{{\emph{i.e.}}\ }
\newcommand{\cf}{{\it cf.}\ }
\newcommand{\eg}{{\it e.g.}\ }
\newcommand{\resp}{{resp.}\ }
\newcommand{\wrt}{{with respect to}}
\newcommand{\lhs}{{left-hand side}}
\newcommand{\rhs}{{right-hand side}}

\newcommand{\dst}{\displaystyle}

\newcommand{\C}{\mathbb{C}}
\newcommand{\HH}{\mathbb{H}}

\newcommand{\Q}{\mathbb{Q}}
\newcommand{\R}{\mathbb{R}}

\newcommand{\Z}{\mathbb{Z}}

\newcommand{\cB}{\mathcal{B}}

\newcommand{\cP}{\mathcal{P}}

\newcommand{\pp}[1]{^{[#1]}}

\DeclarePairedDelimiter\abs{\lvert}{\rvert}%

\DeclarePairedDelimiter\mean{\langle}{\rangle}

\newcommand{\Res}{\operatorname{Res}\!}

\newcommand{\defeq}{\coloneqq} 

\newcommand{\col}{\colon\thinspace}          

\newcommand{\gB}{\mathscr B}       
\newcommand{\gN}{\mathscr N}       
\newcommand{\gQ}{\mathscr Q}       
\newcommand{\gS}{\mathscr S}       
\newcommand{\gL}{\mathscr L}       

\newcommand{\eith}{\ee^{\I\tht}}

\newcommand{\begla}{\begin{equation}}
\newcommand{\beglab}[1]{\begin{equation}	\label{#1}}
\newcommand{\edla}{\end{equation}}

\newcommand{\DD}[2]{{\De\hspace{-.65em}\raisebox{.28ex}{$\scriptscriptstyle /$}\hspace{.25em}^{#1}_{#2}}}

\newcommand{\iimp}{\;\Rightarrow\;}
\newcommand{\imp}{\ens\Rightarrow\ens}
\newcommand{\Imp}{\quad\Rightarrow\quad}

\newcommand{\ti}{\tilde}

\newcommand{\ens}{\enspace}

\newcommand{\bull}{\raisebox{.3ex}{\scalebox{.7}{$\bullet$}\!}}

\newcommand{\ttpd}{{\tfrac{\pi}{2}}}
\newcommand{\tpd}{{\frac{\pi}{2}}}
\newcommand{\trpd}{{\frac{3\pi}{2}}}
\newcommand{\tpq}{\frac{\pi}{4}}
\newcommand{\Itpd}{\frac{\I\pi}{2}}
\newcommand{\Itpq}{\frac{\I\pi}{4}}
\newcommand{\Ittpq}{\frac{3\I\pi}{4}}
\newcommand{\trpq}{\frac{3\pi}{4}}

\newcommand{\Rp}{\R_{>0}}
\newcommand{\Rnn}{\R_{\ge0}}
\newcommand{\Zp}{\Z_{>0}}
\newcommand{\Znn}{\Z_{\ge0}}

\newcommand{\NTL}{^{\operatorname{nt}\!}}
\newcommand{\tnt}{\xrightarrow{\operatorname{n.t.\,}}0}
\newcommand{\tnta}{\xrightarrow{\operatorname{n.t.\,}}\al}
\newcommand{\falM}{f_{\!\frac{\al}{M}}}
\newcommand{\UfalM}{(U_M f)_{\!\frac{-1}{\al M}}}
\newcommand{\fuM}{f_{\!\frac{1}{M}}}
\newcommand{\fbM}{f_{\!\frac{b}{M}}}

\newcommand{\bQalM}{\ov\gQ_{f,M}}
\newcommand{\bQalMU}{\ov\gQ_{U_M f,M}}
\newcommand{\QalM}{\gQ_{f,M}}

\newcommand{\ev}{^{\operatorname{ev}}}
\newcommand{\od}{^{\operatorname{od}}}

\newcommand{\lcm}{\operatorname{lcm}}
\newcommand{\denom}{\operatorname{denominator}}
\newcommand{\SU}{\operatorname{SU}}
\newcommand{\SL}{\operatorname{SL}}

\newcommand{\md}{_\text{\scalebox{.8}{\textup{med}}}}

\newcommand{\Kronecker}[2]{\left(\frac{#1}{#2}\right)}

\newcommand{\datestamp}{{\small{File:\;\hbox{\tt\jobname.tex}
\; \DTMnow \,(Paris time)}}}

\title[Resurgence and Partial Theta Series]{Resurgence and Partial Theta Series}

\author[L.~Han, Y.~Li, D.~Sauzin, S.~Sun]{L.~Han$^{\,1,2}$, Y.~Li$^{\,3}$, D.~Sauzin$^{\,4,1}$, S.~Sun$^{\,1,5}$}

\thanks{$^1$ Department of Mathematics, Capital
  Normal University, Beijing 100048, China}
\thanks{$^2$ Yanqi Lake Beijing Institute of Mathematical Sciences and Applications (BIMSA), Beijing, China}
\thanks{$^3$ Chern Institute of Mathematics and LPMC, Nankai
    university, Tianjin 300071, China}
\thanks{$^4$ CNRS -- 
    Observatoire de Paris, PSL Research University, 75014 Paris,
    France}
\thanks{$^5$ Academy for Multidisciplinary
    Studies, Capital Normal University, Beijing 100048, China}

  \thanks{Corresponding author: {\tt David.Sauzin@obspm.fr}}
   \thanks{\datestamp}


\usepackage[backend=bibtex]{biblatex}
\begin{document}

\maketitle

\begin{abstract}

  We consider partial theta series associated with periodic sequences
  of coefficients, of the form
    $\Th(\tau) := \sum_{n>0} n^\nu f(n) \ee^{\I\pi n^2\tau/M}$,
    with $\nu\in\Znn$ and an $M$-periodic function
    $f : \mathbb{Z} \rightarrow \mathbb{C}$.
    Such a function is analytic in the half-plane $\{\IM\tau>0\}$ and
    as~$\tau$ tends non-tangentially to any $\al\in\Q$, a formal power series appears in the asymptotic behaviour of~$\Th(\tau)$, depending on the
    parity of~$\nu$ and~$f$.
    We discuss the summability and resurgence properties of these series by means of explicit
    formulas for their formal Borel transforms, and the
    consequences for the modularity properties of~$\Th$, or its
    ``quantum modularity'' properties in the sense of Zagier's recent
    theory.
    %
    %
    The Discrete Fourier Transform of~$f$ plays an unexpected role and
    leads to a number-theoretic analogue of \'Ecalle's ``Bridge
    Equations''.
    The motto is: (quantum) modularity $=$ Stokes phenomenon $+$ Discrete Fourier Transform.
\end{abstract}

\bg

This text gathers selected results from our
forthcoming paper~\cite{HLSS}. Here, we aim at giving the global picture,
omitting the most technical details of the
proofs for the sake of clarity.
We are interested in analytic functions defined in
$\HH \defeq \{\, \tau\in\C \mid \IM\tau>0 \,\}$ as \emph{partial theta
  series}, \ie of the form
\beglab{eqdefThan}
\tau\in\HH \mapsto \Th(\tau) \defeq \sum_{n\ge1} a_n \, \ee^{\I\pi n^2\tau/M},
\edla
mostly in the case when the coefficients~$a_n$ are of
the form $n^\nu f(n)$ and $f\col\Z\to\C$ is $M$-periodic.
We will show that their boundary behaviour is better understood when
viewed from the perspective of \'Ecalle's \emph{Resurgence Theory}
\cite{Eca81,MS16}.

Resurgence Theory originates in local analytic dynamics (\'Ecalle-Voronin
theory of parabolic points and Stokes phenomena for nonlinear ODEs
\cite{Eca81,DS14,DS15,mouldSN,Ilya}) and has recently gained
recognition in topology and mathematical physics---see \cite{AM} and
\cite{KS} for the most recent examples. It defines certain subspaces of the space of
all formal series in one indeterminate and shows how to endow them
with a Fr\'echet algebra structure \cite{Sau15};
it then relies on a systematic use of the Laplace transform to relate
certain analytic functions and their divergent asymptotic expansions,
while providing us with the apparatus of \emph{Alien calculus} to
analyse Stokes phenomena (linear or nonlinear).

In our case, the Stokes phenomenon will prove to be related to the
action of the Discrete Fourier Transform on~$f$ and, according to the parity
of~$\nu$ and~$f$, to classical modularity or quantum modularity
properties (in the sense of Zagier~\cite{DZ_QM}).
We shall thus recover and deepen the analytic aspects of various results found in the literature
\cite{AM,CG,GO,GMP,LZ,Shimura,DZ_Top}.

\section{From partial theta series to generating series}

Suppose for the moment that $M\in\Rp$,
$(a_n)_{n\ge1}$ is a sequence of complex numbers and there
exists $K>0$ such that $\abs{a_n}\le K^n$ for all $n\ge1$.
Since $\abs{\ee^{\I\pi n^2\tau/M}} = \ee^{-\pi n^2(\IM\tau)/M}$, the
function~$\Th$ defined by~\eqref{eqdefThan} is analytic on~$\HH$ and exponentially small
as $\IM\tau\to+\infty$.
%

\begin{figure}
\centering
\includegraphics[width=0.45\textwidth]{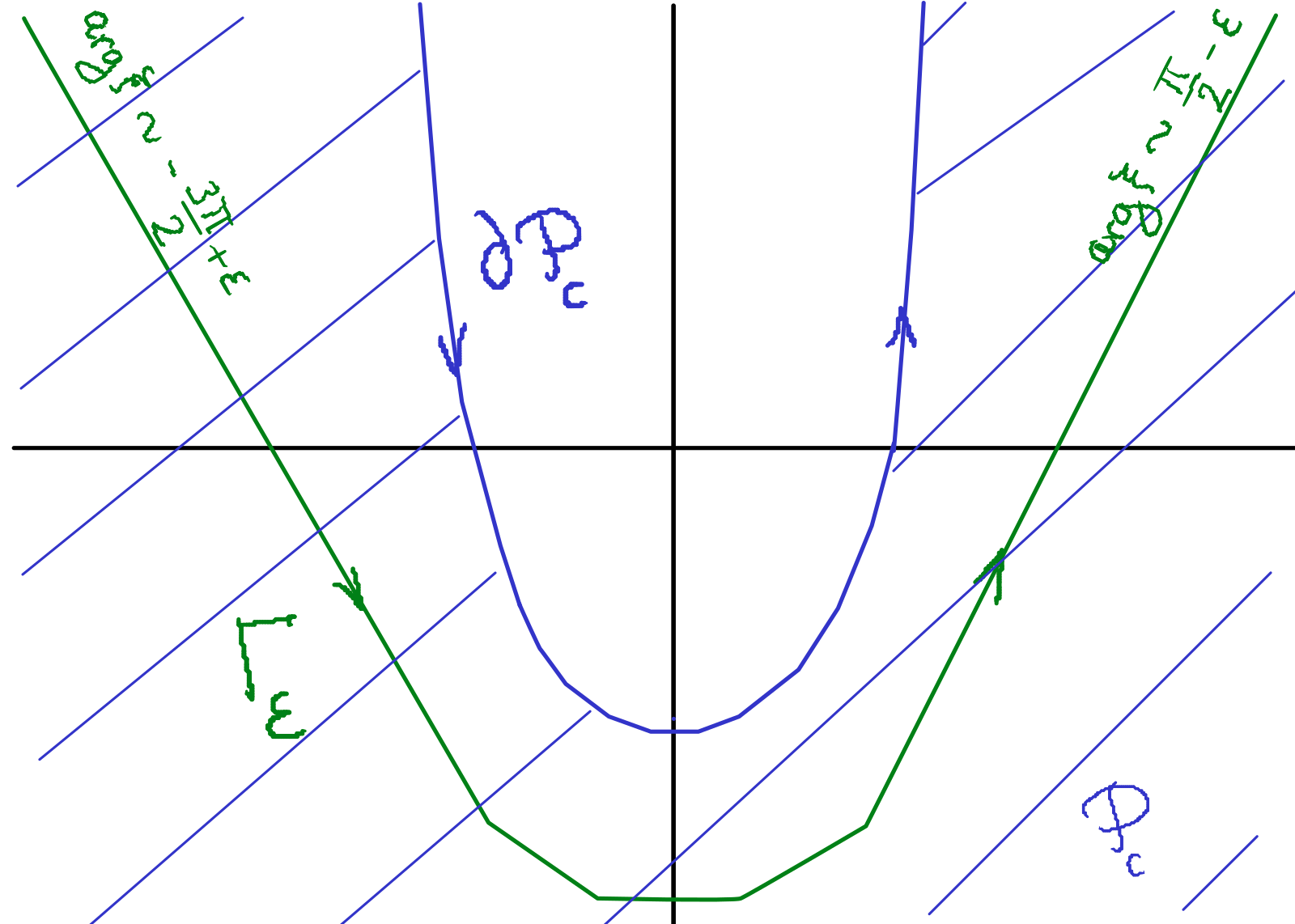}
\caption{The oriented parabola $\pa\cP_c$ and the path~$\Ga_\eps$ in
  the $\xi$-plane.}
\label{fig_parab}
\end{figure}

\begin{thm}    \label{thm_Th_genexp}
Take any $c\in\Rp$ large enough so that the sum of the generating
series
\begla
F(t) \defeq \sum_{n\ge1} a_n \, \ee^{-nt}
\edla
is bounded on $\{\RE t \ge c\}$ and analytic in a neighbourhood of this half-plane.
Then the function
\beglab{eqdefhatphifromF}
\hat\phi(\xi) \defeq \pi^{-1/2} \xi^{-1/2} F(C\,\xi^{1/2})
\qquad \text{with} \ens C \defeq \big(\tfrac{4\pi}{M}\big)^{1/2} \ee^{\I\pi/4},
\edla
where the branch of~$\xi^{1/2}$ is specified by
$-\trpd < \arg\xi < \tpd \iimp -\trpq < \arg(\xi^{1/2}) < \tpq$,
is bounded on $\cP_c \defeq \{\, \RE(C\,\xi^{1/2})\ge c \,\}$
and analytic in a neighbourhood of this set, and we have
  \begla
  \Th(\tau) = \dem \tau^{-1/2} \int_{\pa\cP_c} \ee^{-\xi/\tau} \hat\phi(\xi)
  \,\dd\xi
  \quad \text{for all}\ens \tau\in\HH,
  \edla
  where the parabola $\pa\cP_c$ is oriented anticlockwise---see
  Figure~\ref{fig_parab}---and we choose the branch of the square root
  so that $0<\arg(\tau^{1/2})<\tpd$.
\end{thm}

Notice that when $\xi$ is on the parabola $\pa\cP_c$ and
$\abs\xi\to\infty$, $\arg\xi$ tends to $-\trpd$ or~$\tpd$ so that
$\abs{\ee^{-\xi/\tau}}$ decays exponentially as soon as $0 <
\arg\tau<\pi$.
If $0<\eps<\tpd$ then, for $\arg\tau \in (\eps, \pi-\eps)$, the Cauchy
theorem allows us to replace the integration contour~$\pa\cP_c$ with a
contour $\Ga_\eps\subset \cP_c$ beginning like a half-line of
direction $-\trpd+\eps$ and finishing like a half-line of direction
$\tpd-\eps$, as on Figure~\ref{fig_parab}:
  \beglab{eqThintGaeps}
  \eps < \arg\tau < \pi-\eps \imp
  \Th(\tau) = \dem \tau^{-1/2} \int_{\Ga_\eps} \ee^{-\xi/\tau} \hat\phi(\xi)
  \,\dd\xi.
  \edla
Theorem~\ref{thm_Th_genexp} is implicit in~\cite{GMP} and~\cite{AM}. A
direct proof is obtained from the identity
\begla
  \tau^{1/2} \,\ee^{\sig^2\tau} = \dem \pi^{-1/2} \Big( -\int_0^{\ee^{\I(-\trpd+\eps)}\infty}
  + \int_0^{\ee^{\I(\tpd-\eps)}\infty} \Big)
  \ee^{-\xi/\tau} \, \xi^{-1/2} \,\ee^{-2\sig\xi^{1/2}} \,\dd\xi
\edla
(for any $\sig\in\C^*$), which follows from Borel-Laplace
summation of the \emph{convergent} Puiseux expansion
$\sum \frac{\sig^{2p}}{p!} \tau^{p+\demi}$
of the \lhs\ at~$0$
(indeed: the Borel transform of this series happens to be the odd part
of $\Psi(\xi^{1/2})\defeq \pi^{-1/2} \xi^{-1/2}
\,\ee^{-2\sig\xi^{1/2}}$, we thus recover the
\lhs\ by applying the Laplace transform~\eqref{eqdefLtht} to
$\big[\Psi(\xi^{1/2})-\Psi(-\xi^{1/2})\big]/2$
in any direction, \eg $\tpd\pm\eps$, but we also have
$\gL^{\tpd+\eps}\big[\Psi(-\xi^{1/2})\big] = \gL^{-\trpd+\eps}\big[\Psi(\xi^{1/2})\big]$).
Alternatively, one can use Lemma~1 of~\cite{FN}.

\section{Partial theta series associated with periodic sequences}

From now on, in~\eqref{eqdefThan} we consider sequences $(a_n)$ of the form
$a_n = n^\nu \,f(n)$ with
%
%
%
\beglab{hypodMchi}
\nu\in\Znn, \qquad M\in\Zp, \qquad f \col \Z \to
\C \ens\text{$M$-periodic function.}
\edla
For each $(\nu,f,M)$, the corresponding partial theta series is
\begla
\tau\in\HH \mapsto \Th(\tau; \nu,f,M) \defeq
\sum_{n\ge1} n^\nu \,f(n) \,\ee^{\I\pi n^2\tau/M}
\edla
and we are particularly interested in the local behaviour of this
analytic function as~$\tau$ approaches a point $\al\in\Q$ of the
boundary of~$\HH$.

\begin{rem}   \label{remshiftal}
  For any $\al\in\Q$, we can write
\beglab{eqdeffalM}
\Th(\al+\tau; \nu,f,M) =
\Th\big( \tfrac{M_\al}{M}\tau; \nu,\falM,M_\al \big)
\quad
\text{with notation}\ens
f_{\!\be}(n) \defeq f(n) \, \ee^{\I\pi n^2\be} \ens\text{for $\be\in\Q$},
\edla
where $M_\al$ is a period of~$\falM$,
for instance $\lcm\!\big(M,2\denom(\al/M)\big)$. 
Therefore, the local study near $\tau=0$ of
$\Th(\tau;\nu,f,M)$
with arbitrary $M$ and~$f$ will give us access to the local behaviour
of $\tau\mapsto\Th(\tau; \nu,f,M)$ near any rational.
\end{rem}

Instead of $\Th(\tau;\nu,f,M)$, we will often write
$\Th(\tau;\nu,f)$ or $\Th(\tau)$, omitting some of the parameters if they are
clear from the context.
The cases $\nu=0$ and $\nu=1$ will be the most important ones.

\begin{rem}   \label{remshiftal}
  Note that,
for $\mu\ge0$ and $\nu=2\mu$ or $2\mu+1$,
\beglab{eqdecreaseddtau}
\Th(\tau;2\mu,f,M) = \big(\tfrac{M}{\I\pi}\,
\tfrac{\dd\,}{\dd\tau}\!\big)^{\!\mu} \,\Th(\tau;0,f,M),\quad
\Th(\tau;2\mu+1,f,M) = \big(\tfrac{M}{\I\pi}\, \tfrac{\dd\,}{\dd\tau}\!\big)^{\!\mu} \,\Th(\tau;1,f,M).
\edla
Consequently, even if we state some results with $\nu=0$ or $\nu=1$ as
$\tau\to0$ only, they entail \emph{results for arbitrary~$\nu$ with~$\tau$ approaching any
$\al\in\Q$}, simply by putting together~\eqref{eqdeffalM} and~\eqref{eqdecreaseddtau}.
\end{rem}

\subsubsection*{A few examples}

\begin{enumerate}[topsep=2pt, itemsep=3pt, left=0pt, label=--] 

\item The classical Jacobi theta function appears as
  $\tht_3=1+2\Th(\tau;0,1,1)$.

\item The even primitive quadratic Dirichlet character
  $f(n)=\Kronecker{12}{n}$ (Kronecker symbol), listed as
  $f=\href{http://www.lmfdb.org/Character/Dirichlet/12/11}{\chi_{12}(11,\,\cdot\,)}$
  in~\cite{lmfdb}, gives rise to the Dedekind eta-function
  $\eta(\tau)=\Th(\tau;0,f,12)$, and its Eichler integral
  $\ti\eta(\tau)=\Th(\tau;1,f,12)$ has been investigated particularly
  in~\cite{DZ_Top} and~\cite{CG}.
    %
    %

  \item  The Rogers-Ramanujan identity~\cite{DZ_dilog} involves
    $\tht_{5,j}(\tau) = \Th(\tau;0,f_{5,j},20)$, where
    $f_{5,1}$ and 
    $f_{5,2}$ 
    respectively are the real and imaginary part of
    $\href{http://www.lmfdb.org/Character/Dirichlet/20/7}{\chi_{20}(7,\,\cdot\,)}$,
    an even primitive Dirichlet character \cite{lmfdb}.



\item In Chern-Simons theory with gauge group $\SU(2)$ or in the
  context of Witten-Reshetikhin-Turaev invariants for the Poincar\'e
  homology sphere \cite{LZ,GMP}, there is a partition function
  involving $\Th(\tau;0,f_+,60)$ with
  $f_+\defeq\RE\big(\chi_{60}(23,\,\cdot\,)\big)$, where
  $\href{http://www.lmfdb.org/Character/Dirichlet/60/23}{\chi_{60}(23,\,\cdot\,)}$
  is an odd primitive Dirichlet character \cite{lmfdb}.


\item More generally, for any oriented Seifert fibered integral
  homology sphere with $r\ge3$ exceptional fibers, Andersen and
  Misteg{\r{a}}rd show in~\cite{AM} that the ``Gukov-Pei-Putrov-Vafa
  invariant'' is essentially a partial theta series as
  in~\eqref{eqdefThan}, 
  whose limit as $\tau\to 1/k$ is related to the
  level~$k$ Witten-Reshetikhin-Turaev invariant, at least for $k>0$
  even.
  One can check that, in this case, $a_n$ is a finite sum of odd sequences
  of the form $n^\nu f\pp\nu(n)$, with periodic functions $f\pp\nu$ of the same period
  and~$\nu$ ranging from~$0$ to $r-3$.
  \end{enumerate}

\section{A resurgent divergent series that has three canonical sums}
\label{secdivasympt}

Let $\nu,M,f$ be as in~\eqref{hypodMchi}. With a view to applying
Theorem~\ref{thm_Th_genexp}, we observe that the
generating function has an analytic continuation as a rational function of~$\ee^{-t}$:
\[
F(t) = \sum_{n\ge1} n^\nu\,f(n) \, \ee^{-nt} =
\big(\!\!-\tfrac{\dd\,}{\dd t}\big)^{\!\nu} F_0(t),
\ens
F_0(t) = \sum_{n\ge1} f(n) \, \ee^{-nt} =
\frac{1}{1-\ee^{-M t}} \sum_{1\le\ell\le M} f(\ell)\,e^{-\ell t},
\]
hence we can take $c>0$ arbitrarily small
and $F(t)$ extends from $\{\RE t \ge c\}$ to a meromorphic function
on~$\C$, regular on $\C - \frac{2\pi\I}{M}\Z$.
The mean value $\mean f = \frac{1}{M}\sum\limits_{\ell\in\Z/M\Z}f(\ell)$
controls the regularity at $t=0$:
if $\mean f\neq0$, then~$0$ is a simple pole for~$F_0(t)$ with
residue~$\mean f$;
all other possible poles are of the form $t_n\defeq 2\pi\I n/M$, $n\in\Z^*$.
\sm

To compute the function $\hat\phi(\xi)=\pi^{-1/2} \xi^{-1/2} F(C\,\xi^{1/2})$ of~\eqref{eqdefhatphifromF},
we decompose~$F$ into 
\beglab{eqdecompF}
F(t) = \frac{\nu!\mean f}{t^{\nu+1}} + F^+(t) + F^-(t),\qquad
F^\pm(t)\in\C\{t\},
\qquad F^+\ens\text{even}, \quad F^-\ens\text{odd}.
\edla
%
%
Defining two functions~$\hat\phi^+$ and~$\hat\phi^-$ by
\begin{gather}
  \label{eqdefphipm}
F^+(t) = \pi^{1/2}\, \hat\phi^+\Big(\frac{t^2}{C^2}\Big), \qquad
F^-(t) = \pi^{1/2}\frac{t}{C} \,\hat\phi^-\Big(\frac{t^2}{C^2}\Big),
%
%
\shortintertext{we get}
\label{eqhatphidecsingphipm}
\hat\phi(\xi) = \frac{\nu!\mean f}{\pi^{1/2}C^{\nu+1}} \xi^{-\frac{\nu}{2}-1}
+ \xi^{-1/2}\hat\phi^+(\xi) + \hat\phi^-(\xi),
\qquad \hat\phi^\pm(\xi) \in \C\{\xi\}.
\end{gather}
Notice that $\hat\phi^+$ and~$\hat\phi^-$ are meromorphic on~$\C$ and their possible
poles are of the form
\begla
\xi_n \defeq \frac{t_n^2}{C^2} = \frac{\I \pi n^2}{M} \in \I \Rp, \qquad n\in\Zp.
\edla
The Laplace-like formula~\eqref{eqThintGaeps} thus gives rise to a sum
of three terms:
\begin{align} \label{eqdefThpmd}
\Th(\tau;\nu,f,M) &=
                    \boxed{\, \dem\Ga\big(\tfrac{\nu+1}{2}\big)\mean f
                    \big(\tfrac{\pi}{M}\!\cdot\!\tfrac{\tau}{\I}\big)\!^{-\frac{\nu+1}{2}}
%
                      %
%
+ \Th^+(\tau;\nu,f,M) + \Th^-(\tau;\nu,f,M) } 
   \\[1.5ex]
\text{with}\quad
\Th^+ &\defeq \tfrac{\tau^{-1/2}}{2} \int_{\Ga_\eps} \ee^{-\xi/\tau} \,\xi^{-1/2}\hat\phi^+(\xi)
  \,\dd\xi, \quad
\Th^- \defeq \tfrac{\tau^{-1/2}}{2} \int_{\Ga_\eps} \ee^{-\xi/\tau}\, \hat\phi^-(\xi)\,\dd\xi,
\notag \end{align}
where the first term in~\eqref{eqdefThpmd} has been obtained from the
known value of the Hankel-type
Laplace integral of $\xi^{-\frac{\nu}{2}-1}$, which is $2\pi
\I^{\nu+1}\tau^{-\nu/2}/\Ga(1+\nu/2)$, and the Legendre duplication formula.
%
%

Now, the key point is that the Cauchy theorem allows us to raise the integration
contour~$\Ga_\eps$ until it touches~$0$ and to replace it with the two
half-lines $\ee^{\I(\tpd\pm\eps)}\Rnn$ properly oriented,
taking into account the change of branch for $\xi^{-1/2}$
from $\ee^{\I(-\trpd+\eps)}\Rnn$ to $\ee^{\I(\tpd+\eps)}\Rnn$  in the
case of~$\Th^+$,
whence
\begin{enumerate}[topsep=1.5ex, itemsep=2pt, left=1em, label=\bull] 
\item
  $\frac{\Th^+}{\tau^{-1/2}}$ is the average of the Laplace transforms
  of $\xi^{-1/2}\hat\phi^+(\xi)$ in the directions $\tpd\pm\eps$:
  \beglab{Thpasaver}
    %
    %
\Th^+ = \tau^{-1/2} \times \demi\big(    
\gL^{\tpd-\eps}+\gL^{\tpd+\eps}\big)
\big[\xi^{-1/2}\hat\phi^+(\xi)\big],
  \edla
\item
  $\frac{\Th^-}{\tau^{-1/2}}$ is the difference of the Laplace transforms of $\demi\hat\phi^-(\xi)$ in
  the directions $\tpd\pm\eps$:
  \beglab{eqThmasdiff}
    %
\Th^- = \tau^{-1/2} \times \big(    
\gL^{\tpd-\eps}-\gL^{\tpd+\eps}\big)
\big[\dem\hat\phi^-\big],
  \edla
\end{enumerate}
with the notation
\beglab{eqdefLtht}
  \gL^\tht\hat\ph(\tau) = \int_0^{\eith\infty}
\ee^{-\xi/\tau} \hat\ph(\xi)\,\dd\xi
\edla
in the half-plane $\{\,\arg\tau\in(\tht-\ttpd,\tht+\ttpd)\,\}$.

In the case of $\Th^-(\tau)$, since $\hat\phi^-(\xi)$ is regular at~$0$
and meromorphic, we can go on raising the integration contour as long
as we do not hit the first pole on~$\I \Rp$, which results in an
exponentially small contribution as $\tau\to0$.

In the case of $\Th^+(\tau)$, it is the Puiseux expansion at~$0$ of
\beglab{eqdefphp}
\hat\ph^+(\xi)\defeq\xi^{-1/2}\hat\phi^+(\xi) \in \xi^{-1/2}\C\{\xi\}
%
\edla
that governs the asymptotics:
the standard properties of the Laplace transform~$\gL^{\tpd\pm\eps}$ give
the asymptotic behaviour as $\tau\to0$
non-tangentially, \ie uniformly in any sector
$\eps<\arg\tau<\pi-\eps$, which we denote by $\tau\tnt$.
A few computations lead to formulas involving the analytic
continuation in~$s$ of the $L$-function
$L(s,f) \defeq \sum\limits_{n\ge1} f(n) n^{-s}$
for the Taylor coefficients of $\hat\phi^+(\xi)$ (see \cite{HLSS} for the
details).
The upshot is

\begin{thm}    \label{thm_Th_d_asymp}
  \begin{enumerate}[itemsep=2pt, left=0pt, label=\bull]
    \item
      If the mean value $\mean f$ is non-zero, then the dominant term in~\eqref{eqdefThpmd} is
      the first one
%
  %
  and 
  $\abs{\Th(\tau;\nu,f,M)}\text{\scalebox{.9}{$\xrightarrow[\tau\tnt]{}$}}\infty$. 
    \item  
  The term $\Th^+(\tau;\nu,f,M)$ is the \textbf{median sum} in the
  direction~$\tpd$ of the formal series
  \beglab{eqtiThL}
  \ti\Th(\tau;\nu,f,M) \defeq \sum_{p\ge0} \frac{1}{p!}
  L(-2p-\nu,f) \Big(\frac{\pi\I}{M}\Big)^p \tau^p.
  %
  %
  \edla
  The formal series~$\ti\Th$ is \textbf{resurgent} in $1/\tau$ and its
  Borel transform has multivalued analytic continuation in
  $\C-\{\xi_n\}_{n\ge1}$.
  In particular $\Th^+ \sim_1 \ti\Th$ as $\tau\tnt$.
    \item
  The term $\Th^-(\tau;\nu,f,M)$ is exponentially
  small as $\tau\tnt$: it is $O(\ee^{-c \IM(-1/\tau)})$ with any
  $0<c<\frac{\pi}{M}$
  (much more will be said on it in Section~\ref{secStokesPhen}).
  \end{enumerate}
\end{thm}

A few words on the terminology.
The concept of resurgence \cite{Eca81} is based on the Borel transform
\beglab{eqdefgB}
\gB \col \ti\psi(\tau)=\sum_{p\in \gN}\! c_p \tau^p \mapsto
\hat\psi(\xi)=\sum_{p\in \gN}\! c_p\frac{\xi^{p-1}}{\Ga(p)}
\quad\text{for}\; \gN\!= a+b\Znn\ens \text{with}\ens a,b\in\Rp.
\edla
%
%
The series $\ti\psi(\tau)$ is said to be \emph{resurgent} in~$1/\tau$
if its Borel transform
$\hat\psi(\xi)$ has positive radius of convergence and has ``endless
analytic continuation''---see \cite{Eca81}, or Chapters~5 and~6 of~\cite{MS16}.
If $\abs{\hat\psi(\xi)}$ is bounded by an exponential along a
singularity-free sector $\arg\xi\in[\tht_1,\tht_2]$, then its Laplace
transforms $\gL^\tht\hat\psi$ in the directions
$\tht\in[\tht_1,\tht_2]$ are the analytic continuation of one another and can be glued into one function:
\beglab{eqdefgLinterv}
\gL^{[\tht_1,\tht_2]}\hat\psi(\tau) = \int_0^{\eith\infty}
\ee^{-\xi/\tau} \hat\psi(\xi)\,\dd\xi, \quad
\tht\in[\tht_1,\tht_2]\ens\text{such that}\; \arg\tau\in(\tht-\ttpd,\tht+\ttpd),
\edla
thus yielding the Borel sum $\gS^{[\tht_1,\tht_2]}\ti\psi\defeq \gL^{[\tht_1,\tht_2]}\gB\ti\psi$,
that is analytic in a sectorial neighbourhood of $\tau=0$ with
arguments ranging from $\tht_1-\ttpd$ to $\tht_2+\ttpd$, and has
$1$-Gevrey asymptotic expansion~$\ti\psi(\tau)$ in that domain---this
is the meaning of the symbol ``$\sim_1$'' in
Theorem~\ref{thm_Th_d_asymp}.

If there are singularities along $\eith\Rp$ but not in the sectors
$\arg\xi\in[\tht-\eps,\tht)$ and $(\tht,\tht+\eps]$, then one may resort to
\emph{median} summation $\gS^\tht\md\ti\psi=\gL^\tht\md\gB\ti\psi$,
where $\gL^\tht\md$ is a variant of the Laplace operator introduced by
J.~\'Ecalle;
in general, $\gL^\tht\md\hat\psi$ is not the mere arithmetic average
$\demi(\gL^{\tht-\eps}+\gL^{\tht+\eps})\hat\psi$ but something more
elaborate that involves, in a symmetric way, with well-chosen weights,
the various branches of the analytic continuation of~$\hat\psi$ along
$\eith\Rp - \{\text{singular points}\}$;
see \S1.4 of~\cite{Eca93} or \cite{Menous}.
These weights are such that
$\gS^\tht\md\big(\ti\psi_1\ti\psi_2\big) = \big(\gS^\tht\md\ti\psi_1\big)
\big(\gS^\tht\md\ti\psi_2\big)$, which means that
\beglab{eqgLthtmdmorphism}
\gL^\tht\md\big(\hat\psi_1*\hat\psi_2\big) =
\big(\gL^\tht\md\hat\psi_1\big) \big(\gL^\tht\md\hat\psi_2\big),
\edla
where
$\hat\psi_1*\hat\psi_2(\xi) \defeq \gB\big(\ti\psi_1\ti\psi_2\big)$
is the convolution.\footnote{The
convolution product~$*$ is given by the formula
$\hat\psi_1*\hat\psi_2(\xi) = \int_0^\xi
\hat\psi_1(\xi_1) \hat\psi_2(\xi-\xi_1)\,\dd\xi_1$
for $\abs\xi$ small enough and is known to have endless analytic
continuation when both factors have---see \cite{Eca81}, or Chapter~6
of~\cite{MS16}.
Of course, $\gL^{\tht'}\big(\hat\psi_1*\hat\psi_2\big) = 
\big(\gL^{\tht'}\hat\psi_1\big) \big(\gL^{\tht'}\hat\psi_2\big)$
for a singularity-free direction~$\tht'$, and
$\gS^{\tht'}\big(\ti\psi_1\ti\psi_2\big) = \big(\gS^{\tht'}\ti\psi_1\big)
\big(\gS^{\tht'}\ti\psi_2\big)$ in that case.}

A case in which the median sum coincides with the arithmetic average
is when $\hat\psi(\xi)$ is single-valued for
$\arg\xi\in[\tht-\eps,\tht+\eps]$, \eg a meromorphic function, as is
$\hat\ph^+ = \xi^{-1/2}\hat\phi^+(\xi)$.
Another case is when
$\hat\psi = \hat A*\hat\ph$ with $\hat\ph$ meromorphic and $\hat A$
holomorphic in the sector $\arg\xi\in[\tht-\eps,\tht+\eps]$;
indeed, the identity
$\gL^\tht\md\hat\psi = \demi(\gL^{\tht-\eps}+\gL^{\tht+\eps})\hat\psi$
follows from~\eqref{eqgLthtmdmorphism} because
$\gL^{\tht-\eps}\hat A = \gL^{\tht+\eps}\hat A$ and
$\gL^\tht\md\hat\ph = \demi(\gL^{\tht-\eps}+\gL^{\tht+\eps})\hat\ph$
in that case.

What happens with Theorem~\ref{thm_Th_d_asymp} is
that~\eqref{eqdefphp} defines a function
$\hat\ph^+(\xi) \in \xi^{-1/2}\C\{\xi\}$ that is meromorphic with
singularities in direction~$\tpd$ only (in particular it is endlessly
continuable), so~\eqref{Thpasaver} can be rewritten
$\Th^+(\tau) = \tau^{-1/2} \gL^\tpd\md\hat\ph^+$,
%
%
%
which implies that the function
$\psi^+ \defeq \tau\Th^+(\tau)$
can be written $\psi^+(\tau)= \tau^{1/2} \gL^\tpd\md\hat\ph^+$.
Using $\tau^{1/2} = \gL^\tht\hat A$ with $\hat A(\xi) \defeq \frac{\xi^{-1/2}}{\Ga(1/2)}$ for
any~$\tht$ close to~$\tpd$, we get
$\psi^+(\tau)= \gL^\tpd\md\hat\psi^+$
with
$\hat\psi^+(\xi)\defeq \hat A * \hat\ph^+(\xi)
\in \C\{\xi\}$.
Multiplying $\gL^\tpd\md\hat\psi^+$ by~$\tau\ii$, we finally get
\beglab{eqThpasLmedpsip}
\Th^+(\tau) = \tau\ii\psi^+(\tau) =
\hat\psi^+(0) + \gL^\tpd\md\big(\tfrac{\dd\hat\psi^+}{\dd\xi\;\;}\big)(\tau)
\sim_1 \hat\psi^+(0) +
\gB\ii\big(\tfrac{\dd\hat\psi^+}{\dd\xi\;\;}\big)(\tau)
%
\eqqcolon \ti \Th(\tau).
\edla
We end up with
\begla
\Th^+(\tau) = \gL^\tpd\md\hat\Th = \gS^\tpd\md\ti\Th,
\qquad
\hat\Th = \hat\psi^+(0) \,\de + \frac{\dd\hat\psi^+}{\dd\xi\;\;},
\qquad
\hat\psi^+(\xi) = \frac{\xi^{-1/2}}{\Ga(1/2)} * \big[\xi^{-1/2}\hat\phi^+(\xi)\big]
\edla
by extending~\eqref{eqdefgB} and setting $\cB 1 = \de$ (``Dirac mass at~$0$").

In the Borel plane, $\hat\phi^+(\xi)$ and~$\xi^{-1/2}\hat\phi^+(\xi)$ are meromorphic and hence
endlessly continuable.
Now, convolution with $\xi^{-1/2}$ does not preserve meromorphy, but it
preserves endless continuability, and so does
$\tfrac{\dd\;}{\dd\xi}$, therefore $\hat\Th$ is endlessly
continuable and the formal series~$\ti\Th$ is thus proved to be resurgent.
To check
that~$\ti\Th$ defined by~\eqref{eqThpasLmedpsip} is given
by~\eqref{eqtiThL}, one just needs to compute the Taylor expansion at
$\xi=0$
of~$\hat\psi^+$---we omit the details here and refer the reader
to~\cite{HLSS}.



\begin{rem}
Notice that the resurgent formal series $\ti\Th(\tau)$ has,
for $\tau\in\HH$, three ``canonical'' sums:
$\gL^{\tpd+\eps}\hat\Th(\tau)$, $\gL^{\tpd-\eps}\hat\Th(\tau)$ and
$\gL^\tpd\md\hat\Th(\tau)$, the third of which is the function~$\Th^+(\tau)$.
Interesting phenomena are observed when one compares these three functions---see
Section~\ref{secStokesPhen}.
Notice also that if~$f$ is real-valued, then the coefficients
of~$\ti\Th(\tau)$ \wrt~$\tau/\I$ are real; the realness of
$\Th^+(\tau) = \gS^\tpd\md\ti\Th(\tau)$ for $\tau\in\I\Rp$ in that
case can be viewed as a consequence of a general property of
\'Ecalle's median summation
(anyway, the decomposition~\eqref{eqdefThpmd}
respects realness).
\end{rem}

\begin{rem}   \label{remNTL}
  Using Remark~\ref{remshiftal}, we get from
  Theorem~\ref{thm_Th_d_asymp} the asymptotic behaviour of the
  function $\Th(\tau;\nu,f,M)$ as $\tau\tnta\in\Q$:
  \beglab{eqasymptThetaal}
  \Th(\al+\tau;\nu,f,M) =
\dem\Ga\big(\tfrac{\nu+1}{2}\big)\mean{\falM}
                    \big(\tfrac{\pi}{M}\!\cdot\!\tfrac{\tau}{\I}\big)\!^{-\frac{\nu+1}{2}}
%
  %
%
  + \gS^{\tpd}\md\, \ti\Th_\al(\tau)
  + O(\ee^{-c \IM(-1/\tau)}) \quad\text{as}\ens \tau\tnt
  \edla
  %
  %
  %
  with arbitrary $0<c<\frac{\pi M}{M_\al^2}$,
  and with
\begla
\ti\Th_\al(\tau) \defeq \sum_{p\ge0} \frac{1}{p!} L(-2p-\nu,\falM) \Big(\frac{\pi\I}{M}\Big)^p \tau^p.
\edla
We are thus led to define
\begla
\bQalM \defeq \{\, \al\in\Q \mid \mean{\falM} = 0 \,\},
\edla
so that $\al\notin\bQalM \iimp
\abs{\Th(\tau;\nu,f,M)} \to \infty$ as $\tau\tnta$,
\begla
\hspace{-4.5em}
\al\in\bQalM \iimp \Th\NTL(\al;\nu,f,M) \defeq
\lim_{\tau\tnta}\Th(\tau;\nu,f,M) = L(-\nu,\falM).
\edla
Note that~$f$ and~$\falM$ have same parity, in particular
\begla
\text{$f$ odd} \Imp \bQalM = \Q.
\edla
Note also that the full domain of definition of the boundary function
$\Th\NTL$ is $\bQalM\cup\{\I\infty\}$, since $\Th(\tau;\nu,f,M)\to0$
as $\IM\tau\to\infty$.
We will see in the next
section that $L(-\nu,\falM)=0$ if $\nu$ and~$f$ are odd, and also if $\nu\ge2$
and~$f$ are even, while $L(0,\falM)=-f(0)/2$ if~$f$ is even.
  %
%
\end{rem}

\section{The role of parity}

It so happens that the decomposition of~$f$ into even and odd parts,
\begla
f = f\ev+f\od, \quad
f\ev\!(n)\defeq \tfrac{f(n)+f(-n)}{2}, \quad
f\od(n)\defeq \tfrac{f(n)-f(-n)}{2},
\edla
relates to the decomposition~\eqref{eqdecompF}: for $\nu=0$ one finds
$F_0 = F_0\ev + F_0\od$ with
\beglab{eqFzodFzev}
F_0\ev\!(t) = -\dem f\ev\!(0) +
\tfrac{1}{1-\ee^{-M t}} \sum_{\ell=0}^{M-1} f\od(\ell)\,e^{-\ell
  t},
\quad
F_0\od(t) = \tfrac{1}{1-\ee^{-M t}} \sum_{\ell=0}^{M-1} f\ev\!(\ell)\,e^{-\ell t},
\edla
from which one gets $F^\pm$ by applying $\big(\!\!-\frac{\dd\,}{\dd
  t}\big)^\nu$ and removing $\frac{\nu!\mean f}{t^{\nu+1}}$.
Therefore, when $f$ is even or odd, only one of the two
functions $\Th^+$ and~$\Th^-$ can be nonzero, with the sole exception
of \{$\nu=0$, $f$ even\}, in which case $\Th^+(\tau)=-\demi f(0)$. In general,
%
%
\begin{enumerate}[topsep=1ex, itemsep=1ex, left=1em, label=\bull]
\item
  $\ti\Th(\tau;0,f) = -\demi f\ev\!(0) + \sum_{p\ge0} \frac{1}{p!}
  L(-2p,f\od) \big(\frac{\pi\I}{M}\big)^p \tau^p$
  only depends on~$f\od$ and $f\ev\!(0)$
  and gives rise to $\Th^+(\tau;0,f)$ by median summation,
  while $\Th^-(\tau;0,f)$ only depends on~$f\ev$,
\item
$\ti\Th(\tau;1,f) = \sum_{p\ge0} \frac{1}{p!}
  L(-2p-1,f\ev) \big(\frac{\pi\I}{M}\big)^p \tau^p$ only depends
  on~$f\ev$ and its median sum is
  $\Th^+(\tau;1,f)$, while $\Th^-(\tau;1,f)$ only
  depends on~$f\od$,
%
%
\end{enumerate}
and so on.

\begin{rem}
  If $\nu$ is even (\resp odd), unless~$f$ is even (\resp odd), the
  Borel transform of $\ti\Th(\tau;\nu,f)$ has singularities in the
  direction~$\tpd$ and
  $\ti\Th(\tau;\nu,f)$ is \emph{not} Borel summable in that direction,
  median summation is needed.
  Moreover, the median sum coincides with $\Th(\tau;\nu,f) -
  \dem\Ga\big(\tfrac{\nu+1}{2}\big)\mean f
  \big(\tfrac{\pi}{M}\!\cdot\!\tfrac{\tau}{\I}\big)\!^{-\frac{\nu+1}{2}}$ only if~$f$ is odd
(\resp even).
\end{rem}

\subsubsection*{Discrete Fourier Transform}

The DFT operator~$U_M$ is defined on the space of $M$-periodic
functions $\Z\to\C$ by 
\begla
U_M \col f \mapsto \hat f, \quad
\hat f(n) \defeq \frac{1}{\sqrt M}\, \sum_{\ell\modM}f(\ell)
\ee^{-2\pi\I \ell n/M}
\ens\text{for all}\; n\in\Z.
\edla
This is a unitary operator that respects parity:
\begin{align*}
  &\text{$f$ even} \imp
  \hat f= U_Mf \ens\text{even},\ens
  U_M^2f = f,\ens
  U_M\big(\tfrac{f+\hat f}{2}\big) = \tfrac{f+\hat f}{2}, \ens
  U_M\big(\tfrac{f-\hat f}{2}\big) = -\tfrac{f-\hat f}{2} \\[1ex]
  &\text{$f$ odd} \imp
  \hat f= U_Mf \ens\text{odd},\ens
  U_M^2f = -f,\ens
    U_M\big(\tfrac{f+\I\hat f}{2}\big) = -\I\tfrac{f+\I\hat f}{2}, \ens
    U_M\big(\tfrac{f-\I\hat f}{2}\big) = \I\tfrac{f-\I\hat f}{2}.
\end{align*}
All eigenvectors are of one of the four previous forms, with eigenvalue
$\pm1$ or~$\pm\I$.
\sm

It is~\eqref{eqFzodFzev} that leads us to consider the DFT, since for
all $n\in\Z$,
\beglab{eqresFzevFzod}
\Res\big( F_0\ev(t), t=\tfrac{2\pi\I n}{M} \big) = M^{-\demi} \hat f\od(n), \quad
\Res\big( F_0\od(t), t=\tfrac{2\pi\I n}{M} \big) = M^{-\demi} \hat f\ev(n).
\edla

\section{Number-theoretic Stokes phenomena}   \label{secStokesPhen}

According to Section~\ref{secdivasympt}, the third term of the
decomposition~\eqref{eqdefThpmd}
is the difference of two Laplace transforms:
\begin{multline}   \label{eqThmasStPhen}
\Th^-(\tau;\nu,f,M) =
\tau^{-\demi}\big(\gL^{\tpd-\eps}-\gL^{\tpd+\eps}\big)
\big[\tfrac{\hat\phi^-(\xi)}{2}\big]
= \big(\gL^{\tpd-\eps}-\gL^{\tpd+\eps}\big) \big[\tfrac{\dd\hat\psi^-}{\dd\xi\;\;}\big]\\[1ex]
\quad\text{with} \ens
\hat\psi^-(\xi) \defeq 
\tfrac{\xi^{-1/2}}{\Ga(1/2)} * \tfrac{\hat\phi^-(\xi)}{2}
\in\xi^{\demi}\C\{\xi\}.
\end{multline}
Such a difference usually goes under the name of ``Stokes
phenomenon''.
In the present case, \eqref{eqresFzevFzod} gives rise to
particularly nice formulas, which we spell out only when $\nu=0$ or~$1$
for the sake of clarity: 

\begin{thm}  \label{thmThLmdandStokes}
  We have
  %
%
\begin{align}
  \label{eqdzero}
\Th(\tau;0,f) &=
\dem \hat f\ev\!(0)\big(\tfrac{\tau}{\I}\big)\!^{-1/2}
-\dem f\ev\!(0) + \gS^{\tpd}\md\, \ti\Th(\tau;0,f\od)
+ \big(\tfrac{\tau}{\I}\big)\!^{-\demi}\,\Th(-\tau\ii;0,\hat f\ev) \\
  \label{eqdone}
\Th(\tau;1,f) &=
-\tfrac{M^{1/2}}{2\pi\I} \hat f\ev\!(0) \tau\ii
+ \gS^{\tpd}\md \,\ti\Th(\tau;1,f\ev)
+ \I\big(\tfrac{\tau}{\I}\big)\!^{-\tdemi}\,\Th(-\tau\ii;1,\hat f\od).
\end{align}
\end{thm}

The proof consists in analysing the singularity at $\xi_n=\frac{\I\pi
  n^2}{M}$ of the function~$\hat\phi^-$ of~\eqref{eqdefphipm} 
%
%
so as to write~$\Th^-$
as the sum over $n\ge1$ of one exponentially small contribution for
each singularity.
From~\eqref{eqresFzevFzod} we get
\beglab{eqsingphim}
\dem\hat\phi^-(\xi) = \tfrac{\ee^{\Itpq}\hat
  f\ev(n)}{2\pi\I(\xi-\xi_n)} +\text{reg\quad if $\nu=0$,}
\qquad
\dem\hat\phi^-(\xi) = \tfrac{\I\,\ee^{\Ittpq} n \hat f\od(n)}{-2\pi\I(\xi-\xi_n)^2} +\text{reg\quad if $\nu=1$,}
\edla
where ``$\text{reg}$'' denotes a function that is regular at~$\xi_n$.
The last term in~\eqref{eqdzero}/\eqref{eqdone} thus results from
\beglab{eqexpsmalln}
\big(\gL^{\tpd-\eps}-\gL^{\tpd+\eps}\big) \big[\tfrac{1}{2\pi\I(\xi-\xi_n)}\big] =
\ee^{-\xi_n/\tau}, \quad
\big(\gL^{\tpd-\eps}-\gL^{\tpd+\eps}\big) \big[\tfrac{-1}{2\pi\I(\xi-\xi_n)^2}\big] =
\tau\ii \ee^{-\xi_n/\tau}.
\edla
The other terms in~\eqref{eqdzero}/\eqref{eqdone} are just a rewriting
of the first two terms of~\eqref{eqdefThpmd}, using
$\mean f=\mean{f\ev}=\frac{1}{M^{1/2}}\hat f\ev\!(0)$ and incorporating
the parity information on~$\Th^+$.
By~\eqref{eqdecreaseddtau}, differentiation \wrt~$\tau$
yields a statement for
$\Th(\tau;\nu,f)$
with $\nu\ge2$ as well.
\medskip

We now arrive at the first instance of our motto ``(quantum) modularity $=$ Stokes
phenomenon $+$ DFT'' (according to parity):

\subsection*{Toward Modularity}

The median sum term of~\eqref{eqdzero}, \resp \eqref{eqdone}, is absent
if~$f$ is even, \resp odd.
If $f$ is even, we can rewrite~\eqref{eqdzero} as a weight $\demi$ modularity relation
for a full theta series as far as the modular transformation under
consideration is the negative inversion $S\col \tau\mapsto -\tau\ii$:
\beglab{eqmodulev}
\text{$f$ even},\;\; \tht(\tau;f) \defeq f(0) + 2 \Th(\tau;0,f) \!\imp\!
\tht(\tau;f) =
\big(\tfrac{\tau}{\I}\big)\!^{-\demi}\,\tht(-\tau\ii;\hat f),
\edla
while if $f$ is odd (whence $\hat f\ev(0)=0$), we can
rewrite~\eqref{eqdone} as a weight $\tdemi$ modularity relation:
\beglab{eqmodulod}
\text{$f$ odd} \Imp
\Th(\tau;1,f) = \I\big(\tfrac{\tau}{\I}\big)\!^{-\tdemi}\,\Th(-\tau\ii;1,\hat f)
\quad\text{(exp.\ small as $\tau\tnt$),}
\edla
with a caveat: a priori we have different partial theta series on the
two sides of the relation, one associated with~$f$ and
the other associated with~$\hat f$.
However, if $f$ is an eigenvector of~$U_M$, we do end up with modularity relations
with half-integral weight.
%
%
The classical Jacobi function~$\tht_3$ and the Dedekind $\eta$ function
fall in that category.

If $f$ is a primitive Dirichlet character modulo~$M$,
then $f(-1)= (-1)^\nu$ with $\nu=0$ or~$1$, and $f$ is even ($\nu=0$) or odd
($\nu=1$).
It is then known that $\hat f = \hat f(1) \ov f$, and
\eqref{eqmodulev}--\eqref{eqmodulod} are related to
 the famous results by Shimura on modular forms of half-integral weight~\cite{Shimura}.

Notice that if~$f$ is even or odd but not an eigenvector of~$U_M$, we
still can write
\beglab{eqdefFJ}
  F\defeq\begin{psmallmatrix}f\\ \hat f\end{psmallmatrix},
  \quad
  U_M F \defeq \begin{psmallmatrix} U_M f\\ U_M\hat f\end{psmallmatrix} = J F, \qquad
  J\defeq\begin{psmallmatrix} 0 &1 \\ (-1)^\nu & 0\end{psmallmatrix},
  \quad
  \nu \defeq \text{parity of~$f$},
  \edla
which allows one to construct $2$-dimensional vector-valued modular
forms of half-integral weight:
\[
  \tht(\tau;F) = \begin{psmallmatrix} \tht(\tau;f) \\ \tht(\tau;\hat f)\end{psmallmatrix},
  \qquad
  \Th(\tau;1,F) = \begin{psmallmatrix} \Th(\tau;1,f) \\
    \Th(\tau;1,\hat f) \end{psmallmatrix},
  \]
$\text{$f$ even} \!\imp\!
\tht(\tau;F) =
\big(\tfrac{\tau}{\I}\big)\!^{-\demi}\,J\,\tht(-\tau\ii;F), \;\;
\text{$f$ odd} \!\imp\!
\Th(\tau;1,F) = \I\big(\tfrac{\tau}{\I}\big)\!^{-\tdemi}\,J\,\Th(-\tau\ii;1,F)$.
\medskip

With $\tht_{5,1}(\tau)$ or $\tht_{5,2}(\tau)$, since the real
$20$-periodic functions~$f_{5,j}$ satisfy
$f_{5,1}+\I f_{5,2} = \chi$ with $\chi=\chi_{20}(7,\,\cdot\,)$ even primitive Dirichlet
character, it is better to use $F^*\defeq\begin{psmallmatrix}f_{5,1}\\
  f_{5,2}\end{psmallmatrix}$ instead of the previous~$F$: one finds
$U_{20}\,\chi=\lambda\bar\chi$, $\lambda=\big(\frac{-3 +
  4\I}{5}\big)^{1/4} = s_2+\I s_1$ with
$s_k\defeq \sqrt{\frac{5+(-1)^k\sqrt 5}{10}} = \frac{2}{\sqrt5}\sin\frac{k\pi}{5}$,
hence $U_{20} F^* = S F^*$
with reflection matrix
$S \defeq \begin{psmallmatrix} s_2 & s_1\\
  s_1 & -s_2 \end{psmallmatrix}$.
%
%
Since $f_{5,1}(0)=f_{5,2}(0)=0$, \eqref{eqmodulev} yields
\[
\begin{psmallmatrix} \tht_{5,1}(\tau) \\ \tht_{5,2}(\tau) \end{psmallmatrix} =
\big(\tfrac{\tau}{\I}\big)\!^{-\demi}\, S \,
\begin{psmallmatrix} \tht_{5,1}(-\tau\ii) \\
  \tht_{5,2}(-\tau\ii) \end{psmallmatrix}.
\]

\begin{rem}   \label{remGaussQuadrRecipr}
Let $\al\in\Q^*$. For $f$~even, replacing $\tau$ by $\al+\tau$
in~\eqref{eqmodulev}, we can use~\eqref{eqasymptThetaal} and equate
the leading order terms in the asymptotic behaviour of both sides.
We get
\begla
\mean{\falM} = (\I\al)^{1/2} \mean{\UfalM},
\edla
which stays true even if~$f$ is not even.
This is an identity on quadratic Gauss sums, since
\begla
\mean{\falM} = \frac{1}{M_\al} \,\sum_{n \modMal} f(n) \ee^{\I\pi
  n^2\al/M},
\qquad
\mean{\UfalM} = \frac{1}{M_{\frac{-1}{\al}}} \,\sum_{n \modMmual} (U_M f)(n) \ee^{-\I\pi
  n^2/\al M},
\edla
which is related to the quadratic reciprocity law.
In particular, we see that the
negative inversion $\al\mapsto -1/\al$ maps $\bQalM\cap\Q^*$ to $\bQalMU\cap\Q^*$.
  \end{rem}

\subsection*{Toward Quantum Modularity}

In~\eqref{eqdefThpmd}, the median sum term $\Th^+(\tau;\nu,f,M)$ itself has an interesting
Stokes phenomenon.
In view of \eqref{eqdzero}--\eqref{eqdone}, we restrict our
attention to \{$\nu=0$ and $f$ odd\} or \{$\nu=1$ and $f$ even\},
so $\Th^-\equiv0$, and we consider the exponentially small difference
\begin{gather}
D(\tau) \defeq \big(\gS^{\tpd-\eps}-\gS^{\tpd+\eps}\big) \, \ti\Th(\tau;\nu,f)
= \tau^{-\demi} \big(\gL^{\tpd-\eps}-\gL^{\tpd+\eps}\big) \hat\ph^+
\end{gather}
%
%
%
with $\hat\ph^+ = \xi^{-\demi}\hat\phi^+(\xi)$ (\cf Section~\ref{secdivasympt}),
so that
\beglab{eqLmedStokes}
\Th^+ = \gL^{\tpd}\md\gB\, \ti\Th = 
\gL^{\tpd-\eps}\gB\, \ti\Th(\tau;\nu,f)- \frac{D(\tau)}{2} =
\gL^{\tpd+\eps}\gB\, \ti\Th(\tau;\nu,f)+ \frac{D(\tau)}{2}.
\edla
A residue computation similar to the one that yielded
Theorem~\ref{thmThLmdandStokes} gives
\begin{gather}
  \label{eqsingphp}
\hat\ph^+(\xi) = \tfrac{2\ee^{\Itpq}\hat f\od(n)}{2\pi\I(\xi-\xi_n)} +\text{reg\quad if $\nu=0$,}
\qquad
\hat\ph^+(\xi) = \tfrac{2\I\,\ee^{\Ittpq} n \hat f\ev(n)}{-2\pi\I(\xi-\xi_n)^2} +\text{reg\quad if $\nu=1$,}
\intertext{whence, again by~\eqref{eqexpsmalln},}
\begin{aligned}
  \text{\{$\nu=0$ and $f$ odd\}} &\Imp D(\tau) =
  2 \big(\tfrac{\tau}{\I}\big)\!^{-\demi}\,\Th(-\tau\ii;0,\hat f)
  \\[1ex]
  \text{\{$\nu=1$ and $f$ even\}} & \Imp D(\tau) =
  2 \I\big(\tfrac{\tau}{\I}\big)\!^{-\tdemi}\,\Th(-\tau\ii;1,\hat f).
\end{aligned}
\end{gather}
Therefore, in view of~\eqref{eqLmedStokes}, we can rephrase
\eqref{eqdzero}--\eqref{eqdone} as 
\begin{align}
  \text{$f$ odd} &\imp 
  \Th(\tau;0,f)
  \pm \big(\tfrac{\tau}{\I}\big)\!^{-\demi}\,\Th(-\tau\ii;0,\hat f)
  = \gS^{\tpd\mp\eps}\, \ti\Th(\tau;0,f)
  \\[1ex]
  \text{$f$ even} & \imp  \label{impfevqm}
  \Th(\tau;1,f)
 \pm \I\big(\tfrac{\tau}{\I}\big)\!^{-\tdemi}\,\Th(-\tau\ii;1,\hat f)
  = -\tfrac{M\mean f}{2\pi\I} \tau\ii + \gS^{\tpd\mp\eps}\, \ti\Th(\tau;1,f),
\end{align}
where the \rhs s involve the Borel-Laplace sums of~$\ti\Th(\tau;\nu,f)$
in directions close to~$\tpd$
instead of the median sum in the direction~$\tpd$.
We thus end up with
\begin{thm}   \label{thmModularObstr}
  Suppose \{$\nu=0$, $f$ odd\} or \{$\nu=1$, $f$ even\}.
  If~$f$ is an eigenvector of~$U_M$, put $F\defeq f$ and
  $J\defeq \frac{\hat f}{f} = \pm\I^{\nu+1}$, else
  $F\defeq \begin{psmallmatrix}f\\ \hat f\end{psmallmatrix}$ and
  $J\defeq\begin{psmallmatrix} 0 & 1 \\ (-1)^{\nu+1} & 0\end{psmallmatrix}$.
  Then the ``modular obstruction'' 
  \beglab{eqdefGpm}
  G_\pm(\tau) \defeq
  \Th(\tau;\nu,F)
  \pm
  \I^\nu\big(\tfrac{\tau}{\I}\big)\!^{-\demi-\nu}\,J\,\Th(-\tau\ii;\nu,F)
  +\tfrac{M\mean F}{2\pi\I \tau}
  \edla
  can be written $\gS^{\tht}\,
  \ti\Th(\tau;\nu,F)$ with $\tht=\tpd-\eps$ for~$G_+$ and $\tht=\tpd+\eps$ for~$G_-$,
  and thus extends analytically (by varying~$\tht$ as
  in~\eqref{eqdefgLinterv}) from~$\HH$ through the real
  semi-axis~$\Rp$ to the sector $-2\pi<\arg\tau<\pi$ for~$G_+$, and
  through $\ee^{\I\pi}\Rp$ to the sector $0<\arg\tau<3\pi$ for~$G_-$,
  with $1$-Gevrey asymptotic~$\ti\Th(\tau;\nu,F)$.
\end{thm}

This unexpected regularity property of the modular obstruction~$G_\pm$
was first observed in the case of the Eichler integral~$\ti\eta$ by
D.~Zagier and in the aforementioned case of $\Th(\tau;0,f_+,60)$ in~\cite{LZ},
and gave rise to his theory of quantum modular
forms~\cite{DZ_QM}.
Theorem~\ref{thmModularObstr} is related to the results
of~\cite{GO}---see Section~\ref{sec:new}.

Note the ``built-in modularity'': $G_\pm(\tau)=\pm
\big(\tfrac{\tau}{\I}\big)\!^{-\demi}\,J\,G_\mp(-\tau\ii)$ if $\nu=0$,
and
$G_\pm(\tau)=\pm\I\big(\tfrac{\tau}{\I}\big)\!^{-\tdemi}\,J\,G_\mp(-\tau\ii)$
if $\nu=1$ and $\mean F=0$.
In particular, out of the asymptotics of $G_\pm(\tau)$ as $\tau\to0$ we get
asymptotic expansions for $G_\pm(\tau)$ as $\tau\to\infty$ in terms of
$\ti\Th(-\tau\ii;\nu,F)$.

A function of interest is the non-tangential limit
$\al\in\QalM \mapsto \Th\NTL(\al;\nu,F,M)$, well-defined (according to
Remark~\ref{remNTL}) on
$\QalM\defeq \bQalM \cap \bQalMU$,
%
%
which is $2M$-periodic.
In particular, for~$k$ integer, $\Th\NTL(k;\nu,F,M)$ depends on $k$
mod~$2M$ only and, by taking non-tangential limits
in~\eqref{eqdefGpm}, we get

\begin{cor}
  If $\al\in\Q^*$, then $\al\in\QalM \Leftrightarrow -\al\ii\in\QalM$.
 For $\al=\pm\frac{1}{k}$ with $k\in\QalM\cap\Zp$,
  %
  %
  \begin{align}
    \label{eqNTLod}
  \text{$f$ odd} &\Imp 
  \Th\NTL\big(\!\pm\tfrac{1}{k};0,F\big) =
  -\I\,\ee^{\mp\Itpq} \, k^{1/2} J\,\Th\NTL(\mp k;0,F)
  +G_\pm\big(\!\pm\tfrac{1}{k}\big),  
  \\[1ex]
    \label{eqNTLev}
  \text{$f$ even} & \Imp 
  \Th\NTL\big(\!\pm\tfrac{1}{k};1,F\big) =
  \ee^{\pm\Itpq} \, k^{3/2} \,J\,\Th\NTL(\mp k;1,F)
  \mp \tfrac{M\mean F}{2\pi\I}\, k
  +G_\pm\big(\!\pm\tfrac{1}{k}\big),  
  \end{align}
  %
  %
  with $\dst G_\pm\big(\!\pm\tfrac{1}{k}\big) \sim_1
  \ti\Th\big(\!\pm\tfrac{1}{k};\nu,F\big) =
  \sum_{p\ge0} (\mp1)^p \frac{L(-2p-\nu,F)}{p!} \Big(\frac{\pi\I}{M}\Big)^{\!p} \frac{1}{k^p}$ as $k\to+\infty$.
  \end{cor}
So, the divergent resurgent series~$\ti\Th$ reappears in
the asymptotics of $\Th\NTL\big(\!\pm\tfrac{1}{k};\nu,F\big)$.
Note that in general the dominant terms
in~\eqref{eqNTLod}--\eqref{eqNTLev} are
the first terms of the \rhs s, which are of the form
$k^{\demi+\nu}\times \text{\{$2M$-periodic function of~$k$\}}$.
There are similar formulas for $\Th\NTL(\al\pm\frac{1}{k};\nu,F)$.
See Section~\ref{sec:new} for more on the sets $\QalM$.


\med

Being primitive and real, the Dirichlet character $f(n) =
\Kronecker{12}{n} =\chi_{12}(11,n)$ is an
eigenvector of~$U_{12}$; the eigenvalue is~$1$, we thus recover
weight $\demi$ modularity for the Dedekind $\eta$-function ($\nu=0$)
and quantum modularity for its Eichler integral $\ti\eta$ ($\nu=1$); in
the latter case, our results shed a new light on the resurgence properties already investigated
in~\cite{CG}.
For the Chern-Simons partition function of the Poincar\'e
  homology sphere $\Th(\tau;0,f_+,60)$, 
  recalling that $f_+ = \RE\big(\chi_{60}(23,\,\cdot\,)\big)$
  and introducing $f_- = \IM\big(\chi_{60}(23,\,\cdot\,)\big)$,
  we get
  \begla
F^* \defeq \begin{psmallmatrix}f_+\\
  f_-\end{psmallmatrix}
\Imp
U_{60} F^*=S F^*, \quad S\defeq -\I \begin{psmallmatrix} s_1 & s_2\\
  s_2 & -s_1 \end{psmallmatrix}
\edla
with $s_1,s_2$ as before
(because $\Phi\defeq f_++\I f_-$ is an odd primitive Dirichlet
character and one finds $U_{60}\Phi=(s_2-\I s_1)\bar\Phi$),
whence a $2$D 
quantum modular form $\Th(\tau;0,F^*,60)$.

\section{Alien Derivatives and ``Bridge Equation''}

In this work, the two resurgent series
$\ti\Th(\tau) = \hat\psi^+(0) + \gB\ii\big[\tfrac{\dd\hat\psi^+}{\dd\xi\;\;}\big]
\in\C[[\tau]]$
and $\tau\ii\ti\psi^- = \gB\ii \big[\tfrac{\dd\hat\psi^-}{\dd\xi\;\;}\big]
\in\tau^{\demi}\C[[\tau]]$
have appeared in the decomposition~\eqref{eqdefThpmd} of $\Th(\tau;\nu,f,M)$:
\[
  \Th^+ = \gS^{\tpd}\md\,\ti\Th =
  \hat\psi^+(0) + \gL^{\tpd}\md\big[\tfrac{\dd\hat\psi^+}{\dd\xi\;\;}\big],
  \qquad
  \Th^- = \big(\gL^{\tpd+\eps}-\gL^{\tpd-\eps}\big)\big[\tfrac{\dd\hat\psi^-}{\dd\xi\;\;}\big]
\]
(\cf Theorem~\ref{thm_Th_d_asymp}, \eqref{eqThpasLmedpsip}
and~\eqref{eqThmasStPhen}).
This has led to two interesting Stokes phenomena, via singularity
computations in the Borel plane.

\'Ecalle's Resurgence Theory provides us with an efficient tool to
handle such singularity computations:
with each point~$\om$ of the Riemann surface of the logarithm
is associated the so-called ``alien derivation'' $\De_\om$.
This is an operator which maps a resurgent series~$\ti\ph$ to a resurgent
series whose Borel transform somehow measures the singularities of
certain branches of the analytic continuation of $\gB\ti\ph$ at~$\om$;
in particular, $\De_\om\ti\ph=0$ if all branches of $\gB\ti\ph$ are
regular at~$\om$, which is the case if $\ti\ph$ is a convergent power
series.
The reader is referred to \cite{Eca81,Eca93,MS16} for more details.
In the case of meromorphic singularities, the recipe is motivated
by~\eqref{eqexpsmalln} and can be illustrated with $\xi_n=\frac{\pi
  n^2}{M}\ee^{\Itpd}$ for $n\ge1$: the formulas
\begin{align}
\nu=0\qquad &\Rightarrow&
\De_{\xi_n} \ti\ph^+ &= 2\,\ee^{\Itpq}\hat f\od(n)
&
\De_{\xi_n} \ti\phi^- &= 2\,\ee^{\Itpq}\hat f\ev(n)
\\[.5ex]
\nu=1\qquad &\Rightarrow&
\De_{\xi_n} \ti\ph^+ &= 2\I\,\ee^{\Ittpq} n \hat f\ev(n) \tau\ii
  &
\De_{\xi_n} \ti\phi^- &= 2\I\,\ee^{\Ittpq} n \hat f\od(n) \tau\ii
\end{align}
must be understood as a mere rephrasing of~\eqref{eqsingphp} and~\eqref{eqsingphim}.
For the two resurgent series of our problem, since
$\ti\Th = \tau\ii \ti\psi^+ = \tau^{-\demi}\ti\ph^+$
and $\tau\ii\ti\psi^- = \tau^{-\demi} \ti\phi^-/2$,
the operator~$\De_{\xi_n}$ being a \emph{derivation} annihilating
$\tau^{-\demi}$, we thus get
\begin{align}
\label{eqDetiThzero}
\nu=0\qquad &\Rightarrow&
\De_{\xi_n} \ti\Th &= 2\,\ee^{\Itpq}\hat f\od(n) \tau^{-\demi}
&
\De_{\xi_n} (\tau\ii\ti\psi^-) &= \ee^{\Itpq}\hat f\ev(n) \tau^{-\demi}
\\[1ex]
\label{eqDetiThone}
\nu=1\qquad &\Rightarrow&
\De_{\xi_n} \ti\Th &= 2\I\,\ee^{\Ittpq} n \hat f\ev(n)  \tau^{-\tdemi}
  &
\De_{\xi_n} (\tau\ii\ti\psi^-) &= \I\,\ee^{\Ittpq} n \hat f\od(n) \tau^{-\tdemi}
\end{align}
which means that their Borel transforms have singularities
proportional to $(\xi-\xi_n)^{-\tdemi}$ if $\nu=0$ and
$(\xi-\xi_n)^{-\fdemi}$ if $\nu=1$.

The Stokes phenomena that we have observed just reflect the action of the
\emph{directional alien derivation}
\beglab{eqdirectional}
\DD{}{\tpd}= \sum\limits_{\arg\om=\tpd} \ee^{-\om/\tau}\De_{\om}.
\edla
Indeed, for $\nu=0$: 
\begin{align*} 
%
%
\DD{}{\tpd}\ti\Th(\tau;0,f) &= 2 \big(\tfrac{\tau}{\I}\big)\!^{-\demi}\,\Th(-\tau\ii;0,U_M f\od),
&
\DD{}{\tpd}(\tau\ii\ti\psi^-) &= \big(\tfrac{\tau}{\I}\big)\!^{-\demi}\,\Th(-\tau\ii;0,U_M f\ev),
\intertext{and for $\nu=1$:}
%
%
\DD{}{\tpd}\ti\Th(\tau;1,f) &= 2\I\big(\tfrac{\tau}{\I}\big)\!^{-\tdemi}\,\Th(-\tau\ii;1,U_M f\ev),
  &
\DD{}{\tpd}(\tau\ii\ti\psi^-) &= \I\big(\tfrac{\tau}{\I}\big)\!^{-\tdemi}\,\Th(-\tau\ii;1,U_M f\od).
\end{align*}
In the context of dynamical systems, J.~\'Ecalle has observed that the
action of the alien derivations on divergent series of natural
origin often coincide with the action of a classical operator, typically a
differential operator, giving rise to a ``bridge'' between alien
calculus and classical diffential calculus.
Here, the equations just obtained can be called \emph{Bridge
  Equations} for the problem at hand inasmuch as they show that the
action of the directional alien derivation on our two resurgent series
amount to the action of the DFT operator $U_M$ on the even or odd parts
of~$f$, together with the modular transformation $S\col\tau\mapsto-\tau\ii$.


\begin{rem}
  \'Ecalle's \emph{Alien Calculus} gives a formula for the commutator
  of an alien derivation~$\De_\om$ and the natural derivation
  $\tfrac{\dd\,}{\dd\tau}$, which amounts to the fact that each
  homogeneous component $\ee^{-\om/\tau}\De_{\om}$ in~\eqref{eqdirectional} commutes with the
  natural derivation $\tfrac{\dd\,}{\dd\tau}$, namely
  \begla
  \De_\om \tfrac{\dd\,}{\dd\tau} = \big(\tfrac{\dd\,}{\dd\tau} +
  \om\tau^{-2}\big)\De_\om.
  \edla
  We can thus differentiate~\eqref{eqDetiThzero}--\eqref{eqDetiThone}
  \wrt~$\tau$ and, together with~\eqref{eqdecreaseddtau}, this implies
\begin{align*}
\nu=2 \; &\Rightarrow& 
\De_{\xi_n} \ti\Th(\tau;2,f,M) &= -2\, n^2 \hat f\od(n) \big(\tfrac{\tau}{\I}\big)^{-\fdemi}
                                 +\tfrac{M}{\pi} \hat f\od(n) \big(\tfrac{\tau}{\I}\big)^{-\tdemi},
\\[1ex]
\nu=3 \; &\Rightarrow& 
\De_{\xi_n} \ti\Th(\tau;3,f,M) &= 
                                 -2\I\, n^3 \hat f\ev(n) \big(\tfrac{\tau}{\I}\big)^{-\sdemi}
                               + 3\I \tfrac{M}{\pi}\, n \hat f\ev(n)
                                 \big(\tfrac{\tau}{\I}\big)^{-\fdemi},
\\[1ex]
\nu=4 \; &\Rightarrow& 
\De_{\xi_n} \ti\Th(\tau;4,f,M) &= 2\, n^4 \hat f\od(n) \big(\tfrac{\tau}{\I}\big)^{-\ndemi}
                                 -6\tfrac{M}{\pi}\,n^2 \hat f\od(n) \big(\tfrac{\tau}{\I}\big)^{-\sdemi}
                                 +3\big(\tfrac{M}{\pi})^2 \hat f\od(n) \big(\tfrac{\tau}{\I}\big)^{-\fdemi},
\end{align*}
etc.
\end{rem}


\begin{rem}
  The growth as $p\to\infty$ of
  $[\tau^p]\ti\Th= \frac{1}{p!}  L(-2p-\nu,f)
  \big(\frac{\pi\I}{M}\big)^p =
  \frac{\dd^p \hat\psi^+}{\dd\xi^p\;}\!(0)$ is governed by the
  nearest-to-origin singularity of
  $\hat\psi^+=\gB(\tau\ti\Th)\in\C\{\xi\}$; taking
  $n^*\ge1$ minimum such that $\hat f\od(n^*)$ or $\hat
  f\ev(n^*)\neq0$, we thus have, by virtue
  of~\eqref{eqDetiThzero}--\eqref{eqDetiThone},
  $\De_{\xi_{n^*}}(\tau\ti\Th)=S_\nu\tau^{\demi-\nu}$
  with $S_0 = 2\,\ee^{\Itpq}\hat f\od(n^*)$, $S_1 = 2\I\,\ee^{\Ittpq} n^*
  \hat f\ev(n^*)$,
  which means
  \begla
  \hat\psi^+(\xi) = \tfrac{S_0}{2\Ga(\demi)}(\xi-\xi_{n^*})^{-\demi} +
  \hat R_0(\xi)
  \quad \text{or} \ens
  -\tfrac{S_1}{4\Ga(\demi)}(\xi-\xi_{n^*})^{-\tdemi} + \hat R_1(\xi)
  \edla
  with $\hat R_\nu(\xi)\in\C\{\xi\}$ convergent in a disc of radius larger
  than~$\abs{\xi_{n^*}}$.
  The coefficients of the Taylor expansion of
  $(\xi-\xi_{n^*})^{-\demi-\nu}$ thus give an asymptotic equivalent of
  $[\tau^p]\ti\Th$ up to exponential precision.

  Another commonly used variable is $q=\ee^{2\pi\I\tau}$, and the resulting
  asymptotic expansions in powers of $q-1$ as $q\to1$
  (non-tangentially from within the
  unit disc) can be handled as follows:
  put $Q = \frac{q-1}{2\pi\I} = g(\tau)$, so $\tau = g\ii(Q)\in
  Q+Q^2\C\{Q\}$,
  then the new asymptotic expansion is
  $\ti\Th^* = \ti\Th\circ g\ii$, hence it must be resurgent
  too and the \emph{alien chain rule} gives
  $\De_{\xi_n}\ti\Th^*(Q) =
  \big(\ee^{-\xi_n\big(\frac{1}{\tau}-\frac{1}{Q}\big)}\De_{\xi_n}\ti\Th(\tau)\big)_{\mid
      \tau=g\ii(Q)}$,
  thus providing access to the asymptotics
  of $[Q^p]\ti\Th^*$ as $p\to\infty$.
 Applying this with the quadratic character $f(n) = \Kronecker{12}{n}$, one
 recovers the asymptotic of Glaisher's T-numbers, as well as that of
 Stoimenow's numbers, which count regular linearized chord diagrams of
 degree~$p$~\cite{DZ_Top}.
\end{rem} 


\section{Action of the modular group and (quantum) modularity} \label{sec:new}

The standard action of a general modular transformation
$\ga = \begin{psmallmatrix} a & b\\
  c & d \end{psmallmatrix} \in \SL(2,\Z)$ on~$\HH$,
$\ga(\tau) = \frac{a\tau+b}{c\tau+d}$,
induces an action on our partial theta series.

\subsection*{Parabolic Case}
%
The case
$\ga = \begin{psmallmatrix} 1 & b\\
  0 & 1 \end{psmallmatrix}$ is straightforward;
actually, it is a particular case of~\eqref{eqdeffalM}:
\begla
\Th(\tau+b;\nu,f,M) = \Th\big( \tfrac{M'}{M}\tau; \nu,\fbM,M' \big)
\quad
\text{where}\ens
\fbM \defeq \La_M^b f,\quad \La_M(n) \defeq \ee^{\I\pi n^2/M}
\edla
and $M'$ is a period of~$\fbM$. We note that $\La_M(n+M) =
(-1)^M\La_M(n)$, hence
\begla
\text{$M$ or $b$ even} \Imp
\Th(\tau+b;\nu,f,M) = \Th\big( \tau; \nu,\fbM,M \big)
\edla
(because one can then take $M'=M$).
Notice that if the support of~$f$ is such that
\beglab{assumptsupp}
\text{$\exists n_0\in\Z$ such that, $\forall n\in\Z$,} \ens
  f(n)\neq0 \iimp n^2 = n_0^2 \;\moddM,
  \edla
  then $\fuM = \La_M(n_0) f$ and $\fbM = \La_M^b(n_0) f$.
  Notice also that, since $\La_M^{2M}=1$, in all cases
  \begla
  b=0\;\moddM \Imp \fbM=f
  \edla
  (indeed: $\Th(\tau;\nu,f,M)$ is $2M$-periodic in~$\tau$).

\subsection*{Non-Parabolic Case}
%
If
$\ga = \begin{psmallmatrix} a & b\\
  c & d \end{psmallmatrix}$ with $c\neq0$,
then one can assume $c>0$ without loss of generality (replacing~$\ga$ with
$-\ga$ if necessary).
Using what is known from Section~\ref{secStokesPhen} on the action of
$S=\begin{psmallmatrix} 0 & -1\\
  1 & 0 \end{psmallmatrix}$ and writing
\begla
\ga(\tau) = c\ii \big(a + S(c\tau+d)\big),
\edla
one can compute $\Th\big(\ga(\tau);\nu,f,M\big)$. Here is what one finds,
assuming for the sake of simplicity
\beglab{assumptMevetc}
\text{$M$ even,\; $\nu\in\{0,1\}$,\; $f$ even or odd}
\edla
till the end of this section (otherwise some adjustments are
required):

\begin{thm}
  Assume~\eqref{assumptMevetc} and $\ga = \begin{psmallmatrix} a & b\\
  c & d \end{psmallmatrix} \in \SL(2,\Z)$ with $c>0$. 
Then the formula
\beglab{eqdefhg}
h(n) \defeq(cM)^{-1/2} \La_M^{bd}(n) \sum_{r \modM} f(r+dn) g(r)
\ee^{2\pi\I bnr/M},
\quad
g(r) \defeq\sum_{\substack{\ell\modcM \\ \text{s.t.}\;\ell=r\modM}} \La_{cM}^a(\ell)
\edla
defines an $M$-periodic function~$h$ of the same parity as~$f$
and, using the notation~$\tht(\cdot\,;\,\cdot)$ introduced in~\eqref{eqmodulev},
\begin{align}
\text{$f$ even} & \Imp \qquad\ens \tht\big(\ga(\tau);f) =
                  \big(\tfrac{c\tau+d}{\I}\big)^{\demi} \tht(\tau;h) \label{eqgenclassmodev} \\[1ex]
\text{$f$ odd} & \Imp \Th\big(\ga(\tau);1,f,M) =
                 \I \big(\tfrac{c\tau+d}{\I}\big)^{\tdemi}
                 \Th(\tau;1,h,M) \label{eqgenclassmodod} \\[1ex]
\text{$f$ even} & \Imp 
  %
%
\Th(\tau;1,h,M)  \notag
                  \pm \I \big(\tfrac{c\tau+d}{\I}\big)^{-\tdemi}
                  \Th\big(\ga(\tau);1,f,M\big) \\[.5ex]
& \qquad\qquad\qquad\qquad\qquad\quad\;                 = \tfrac{(cM)^{1/2}}{2\pi} f(0) \tfrac{\I}{c\tau+d} +
                  \gS^{\tpd\mp\eps}\,\ti\Th(c\tau+d;1,\La_{cM}^{-d}
    h,cM) \label{eqgenmodobstrev} \\[1ex]
\text{$f$ odd} & \Imp 
  %
%
\Th(\tau;0,h,M)  \notag
                  \mp \big(\tfrac{c\tau+d}{\I}\big)^{-\demi}
                  \Th\big(\ga(\tau);0,f,M\big) \\[.5ex]
                & \qquad\qquad\qquad\qquad\qquad\qquad\qquad
                  = \gS^{\tpd\mp\eps}\,\ti\Th(c\tau+d;0,\La_{cM}^{-d}h,cM).
                  \label{eqgenmodobstrod}
    \end{align}
\end{thm}

Notice that, in~\eqref{eqgenmodobstrev} and~\eqref{eqgenmodobstrod},
the \rhs\ has analytic continuation through the real axis to the right
or to the left of $\ga\ii(\I\infty)=-d/c$ (according to the~$\pm$
sign) by the standard properties of Borel-Laplace summation.

Imposing to~$\ga$ to belong to a certain congruence subgroup and
making use of classical properties of quadratic Gauss sums, one
obtains a relation between~$h$ and~$f$ simple enough to allow us to interpret the \lhs\
of~\eqref{eqgenmodobstrev} and~\eqref{eqgenmodobstrod} as a genuine
modularity obstruction, with appropriate automorphy factor.
We thus get classical modularity from~\eqref{eqgenclassmodev}--\eqref{eqgenclassmodod}  or quantum modularity
from~\eqref{eqgenmodobstrev}--\eqref{eqgenmodobstrod} according to the
parity of~$\nu$ and~$f$.

Specifically, if we assume $\ga\in\Ga_0(2M) \defeq \{\, c = 0 \moddM
\,\}$,
then $g(r)=0$ in~\eqref{eqdefhg} unless $r=0\modM$ and
\begla
g(0) = (cM)^{1/2} \ee^{\I\pi/4} \eps_d\ii \big(\tfrac{2Mc}{\abs{d}}\big)
\quad\text{with}\quad
\eps_d \defeq \begin{cases} 1 & \text{if $d=1 \modfour$} \\
  \I & \text{if $d=3 \modfour$}
\end{cases}
\edla
and $\big(\tfrac{m}{\abs{d}}\big) \defeq$ Jacobi symbol (notice
that~$d$ is necessarily odd since $ad=1\moddM$),
whence
\[
  \ga \in \Ga_0(2M) \Imp
  h(n) = \ee^{\I\pi/4} \,\eps_d\ii\, \big(\tfrac{2Mc}{\abs{d}}\big)
  \La_M^{bd}(n) f(dn).
  \]
Since here we assume~$M$ even, it follows that
\beglab{eqhcaseGa}
\ga\in\Ga(2M) \Imp
h = \ee^{\I\pi/4}\, \big(\tfrac{2Mc}{\abs{d}}\big) f,
\edla
where $\Ga(2M) \defeq \{\, a = d= 1 \moddM \;\text{and}\;
b =  c = 0 \moddM \,\}$.

If we use the larger subgroup $\Ga_1(2M)\defeq \{\,  a = d= 1 \moddM \;\text{and}\;
c = 0 \moddM \,\}$
instead of $\Ga(2M)$, then we still have a relation similar to one of
the results of~\cite{GO} under the support assumption~\eqref{assumptsupp}:
\beglab{eqhcaseGaunSA}
\ga \in \Ga_1(2M) \; \text{and Assumption~\eqref{assumptsupp} holds} \Imp
h = \La_M^b(n_0)\, \ee^{\I\pi/4}\, \big(\tfrac{2Mc}{\abs{d}}\big) f.
\edla
Another option available when~$f$ is a Dirichlet character modulo~$M$ is to use $\Ga_0^0(2M)\defeq \{\,  b=0 \moddM \;\text{and}\;
c = 0 \moddM \,\}$:
\beglab{eqhcaseGazzDC}
\ga \in \Ga_0^0(2M) \; \text{and $f$ Dirichlet character $\!\modM$} \Imp
h = \ee^{\I\pi/4}\, \eps_d\ii\, \big(\tfrac{2Mc}{\abs{d}}\big) f(d) f.
\edla

\emph{The relations \eqref{eqgenclassmodev}--\eqref{eqgenmodobstrod} thus
give rise to explicit transformation formulas for the function
$\Th(\,\cdot\,;\nu,f,M)$ under the action of the congruence
subgroup~$\Ga$ defined by}
\begin{multline} \label{eqdefGa}
\Ga\defeq\Ga(2M) \ens\textit{or, if Assumption~\eqref{assumptsupp} holds,}\;
\Ga \defeq \Ga_1(2M) \\[1ex]
\textit{or, if $f$ is a Dirichlet character  $\!\modM$,}\;
\Ga \defeq \Ga_0^0(2M).
\end{multline}

\subsection*{The Boundary Function $\Th\NTL$}
%
Let us give more details on the cases
\begla
  \text{\{$\nu=0$ and $f$ odd\}}
  \ens\text{or}\ens
  \text{\{$\nu=1$ and $f$ even\}}
  \edla
  (still with~$M$ even). Recall that the boundary function
  $\Th\NTL(\al) = \Th\NTL(\al;\nu,f,M)$ is defined for $\al\in\bQalM\cup\{\I\infty\}$.

If $\ga= \begin{psmallmatrix} 1 & b\\
  0 & 1 \end{psmallmatrix} \in\Ga$ is parabolic, we have
$\Th\big(\ga(\tau);\nu,f,M\big) = \la(\ga) \Th(\tau;\nu,f,M)$ with
$\la(\ga) \defeq 1$ or $\la(\ga)\defeq \La_M^b(n_0)$ according as which
option holds in~\eqref{eqdefGa}; in particular, $\bQalM$ is stable under~$\ga$.

If $\ga = \begin{psmallmatrix} a & b\\
  c & d \end{psmallmatrix} \in \Ga$ is non-parabolic, then we assume $c>0$ and
\eqref{eqhcaseGa}--\eqref{eqhcaseGazzDC} say that the function~$h$
of~\eqref{eqdefhg} is $h =\la(\ga) f$ with
$\la(\ga) \defeq \ee^{\I\pi/4}\, \big(\tfrac{2Mc}{\abs{d}}\big)$ or
$\la(\ga)\defeq \La_M^b(n_0)\,\ee^{\I\pi/4}\,
\big(\tfrac{2Mc}{\abs{d}}\big)$ or
$\la(\ga)\defeq \ee^{\I\pi/4}\, \eps_d\ii\,
\big(\tfrac{2Mc}{\abs{d}}\big) f(d)$.
      %
      %
The relations~\eqref{eqgenmodobstrev}--\eqref{eqgenmodobstrod} thus
yield
\begin{multline} \label{eqrelqmf}
\Th(\tau;\nu,f,M) 
                  \mp \tfrac{1}{\I\la(\ga)} \big(\tfrac{c\tau+d}{\I}\big)^{-\nu-\demi}
                  \Th\big(\ga(\tau);\nu,f,M\big) \\[.5ex]
                = \tfrac{(cM)^{1/2}}{2\pi\la(\ga)} f(0) \tfrac{\I}{c\tau+d} +
                  \gS^{\tpd\mp\eps}\,\ti\Th(c\tau+d;\nu,\La_{cM}^{-d}f,cM).
\end{multline}
Assuming $f(0)=0$ and taking the limit as~$\tau$ tends
non-tangentially to~$\ga\ii(\I\infty)$ (in which case the second term
tends to~$0$ exponentially fast), or to a rational
$\al>\ga\ii(\I\infty)$ (in which case we use `$-$'), or to a rational
$\al<\ga\ii(\I\infty)$ (in which case we use `$+$'), we get
\begin{thm}
  Suppose $M$ is even and $f(0)=0$.
  If $f$ is even then the domain of definition
  $\bQalM\cup\{\I\infty\}$ of $\Th\NTL$ is stable under~$\Ga$:
  it is a union of cusps, which contains at least the cusp at infinity, and also
  the cusp at~$0$ if $\mean{f}=0$.
If $f$ is odd then $\bQalM=\Q$.
In all cases, the function $\Th\NTL \col \bQalM\cup\{\I\infty\} \to \C$ is a quantum modular
form \wrt~$\Ga$, of weight~$\demi$ if \{$\nu=0$ and $f$ is odd\},
of weight~$\tdemi$ if \{$\nu=1$ and $f$ is even\}.
  \end{thm}

Notice that, as a consequence, $\bQalM$ is dense in~$\Q$ if
$f(0)=0$. 
Indeed, each cusp of~$\Ga$ is dense in~$\Q$
because~$\Ga$ contains $\Ga(2M)$, which is a normal subgroup of finite
index of $\SL(2,\Z)$.

\begin{rem}
  The number
  $\Th\NTL(\al)$ is just the limit value of
  $\Th(\tau;\nu,f,M)$ as
  $\tau\tnta$, but we can also consider the whole asymptotic expansion
  of~\eqref{eqrelqmf} as $\tau\tnta$:
  the formal series $\big(\ti\Theta_\al\big)_{\al\in\bQalM}$
  introduced in Remark~\ref{remNTL} make up
  what is called a \emph{strong quantum modular form} in~\cite{DZ_QM}.
\end{rem}

The case when Assumption~\eqref{assumptsupp} holds has already been
considered (with $M$ even or odd) in~\cite{GO}, and we get the same
congruence subgroup as them in that case.
A notable example where this assumption is fulfilled is,
following~\cite{GM} or~\cite{AM}, the partial theta series giving the Gukov-Pei-Putrov-Vafa
  invariant for an oriented Seifert fibered integral
  homology sphere with~$3$ exceptional fibers (a Brieskorn $3$-sphere).

\begin{rem}
  For arbitrary $\nu\ge2$, since Borel-Laplace (median) summation
  commutes with $\tfrac{\dd\,}{\dd\tau}$, we can
  use~\eqref{eqdecreaseddtau} and derive transformation formulas for $\Th(\tau;\nu,f,M)$,
  not as simple as when $\nu=0$
  or~$1$, but rather pertaining to the theory of \emph{higher depth quantum modular
    forms}~\cite{BKL}.
\end{rem}

\begin{rem}
  \textbf{Examples with $f$ even and $\bQalM=\Q$.}
  Suppose that $\nu=1$ and~$f$ is even.
  If Assumption~\eqref{assumptsupp} holds for~$f$, then we get
  $f_{\frac{\al+1}{M}} = \La_M(n_0) f_{\frac{\al}{M}}$ for every
  $\al\in\Q$, hence
  \beglab{impAsinvT}
  \text{Assumption~\eqref{assumptsupp} holds} \Imp
  \text{$\bQalM$ is invariant under
    $T\col\tau\mapsto\tau+1$.}
  \edla
  
  \noindent Using~\eqref{impfevqm}, we also have
  \begla
  \text{$f$ eigenvector of~$U_M$ and $\mean{f}=0$} \Imp
  \text{$\bQalM\cup\{\I\infty\}$ is invariant under
    $S\col\tau\mapsto-1/\tau$.}
  \edla
  Since $\langle S,T\rangle=\SL(2,\Z)$ acts transitively on $\Q\cup\{\I\infty\}$,
  we deduce that \emph{$\bQalM=\Q$ whenever~$f$ is a zero mean value even
  eigenvector of~$U_M$ such that~\eqref{assumptsupp} holds}.
  This is what happens with the Eichler integral $\ti\eta(\tau) =
  \Th\big(\tau; 1, \chi_{12}(11,\,\cdot\,),12\big)$ of the
  Dedekind~$\eta$ function.
  
  There are cases where~$f$ is not an eigenvector of~$U_M$, but $U_M f$
  is a linear combination of~$f$ and~$\ti f$ with~$f$
  satisfying~\eqref{assumptsupp} with a certain~$n_0$ and~$\ti f$
  satisfying~\eqref{assumptsupp} with a certain~$\ti n_0$, so that
  $\bQalM\cap\bQalMU$ is invariant under~$T$
  by~\eqref{impAsinvT}.
  If moreover $\mean f=f(0)=0$, then~\eqref{impfevqm} shows that
  $\big(\bQalM\cup\{\I\infty\}\big)\cap\big(\bQalMU\cup\{\I\infty\}\big)$
  is invariant under~$S$, whence
  $\bQalM=\bQalMU = \Q$ in thoses cases.
  This is what happens with the Eichler integrals $\Th(\tau;1,f_{5,j},20)$
  of the two modular forms involved in the Rogers-Ramanujan identity.
  \end{rem}

  \vspace{.3cm}


\subsubsection*{Acknowledgements}

The 2nd and 3rd authors thank Capital Normal
  University for its hospitality.
The 2nd author aknowledges support from NSFC (No.11771303).
The 4th author is partially supported by National Key R\&D Program
of China (2020YFA0713300), NSFC (No.s 11771303, 12171327, 11911530092,
11871045).
This paper is partly a result of the ERC-SyG project, Recursive and
Exact New Quantum Theory (ReNewQuantum) which received funding from
the European Research Council (ERC) under the European Union's Horizon
2020 research and innovation programme under grant agreement No
810573.

\vspace{.3cm}


\printbibliography


\end{document}